\documentclass[a4paper, reqno, 11pt]{amsart}

\usepackage{amsfonts, amsmath, amssymb, amsthm, amscd, setspace, subfigure,mathtools,stmaryrd}
\usepackage[heightrounded,textwidth=400pt,top=1in,left=1in,right=1in,bottom=1in]{geometry}
\usepackage{layout}
\usepackage[utf8]{inputenc}
\usepackage[unicode,pagebackref=true]{hyperref}
\usepackage{cite}
\usepackage{tikz}
\usetikzlibrary{calc}

\numberwithin{equation}{section}

\newcommand\be{\begin{eqnarray}}
\newcommand\ee{\end{eqnarray}}

\newcommand{\C}{{\mathbb C}}

\theoremstyle{definition}

\newtheorem{example}{Example}[section]

\author[Bhoja]{Niren Bhoja}
\email{\href{mailto: niren.bhoja@nottingham.ac.uk}{niren.bhoja[at]nottingham.ac.uk}}
\author[Krasnov]{Kirill Krasnov}
\email{\href{mailto: kirill.krasnov@nottingham.ac.uk}{kirill.krasnov[at]nottingham.ac.uk}, ORCID: \href{https://orcid.org/0000-0003-2800-3767}{0000-0003-2800-3767}}
\address{School of Mathematical Sciences, University of Nottingham, Nottingham, NG7 2RD, UK}

\title{Spinors from pure spinors}

\date{December 2023}

\begin{document}

\begin{abstract}\noindent 
We propose and develop a new method to classify orbits of the spin group ${\rm Spin}(2d)$ in the space of its semi-spinors. The idea is to consider spinors as being built as a linear combination of their pure constituents, imposing the constraint that no pair of pure spinor constituents sums up to a pure spinor. We show that this leads to a simple combinatorial problem that has a finite number of solutions in dimensions up to and including fourteen. We call each distinct solution a combinatorial type of an impure spinor. We represent each combinatorial type graphically by a simplex, with vertices corresponding to the pure constituents of a spinor, and edges being labelled by the dimension of the totally null space that is the intersection of the annihilator subspaces of the pure spinors living at the vertices. We call the number of vertices in a simplex the impurity of an impure spinor. In dimensions eight and ten the maximal impurity is two. Dimension twelve is the first dimension where one gets an impurity three spinor, represented by a triangle. In dimension fourteen the generic orbit has impurity four, while the maximal impurity is five. We show that each of our combinatorial types uniquely corresponds to one of the known spinor orbits, thus reproducing the classification of spinors in dimensions up to and including fourteen from simple combinatorics. Our methods continue to work in dimensions sixteen and higher, but the number of the possible distinct combinatorial types grows rather rapidly with the dimension. 

\end{abstract}

\maketitle
\tableofcontents

\section{Introduction} 

The problem of classification of spinors is the problem of finding, in each dimension, all different spinor orbits, and computing the corresponding stabilisers. This problem has two possible versions. One is over $\C$, where one considers the action of ${\rm Spin}(n,\C)$ on the spinor representation $S$. In even dimensions it is natural to consider the action on the irreducible semi-spinor representation $S^+$ instead. For low $n\leq 10$ the solution of this problem has been known for a long time, and will be reviewed below. The classification in dimension twelve was done in \cite{Igusa}. Using similar techniques, the classification of spinors in dimension thirteen was achieved in \cite{Kac-Vinberg}. 
The case of fourteen dimensions was treated in \cite{Popov}, see also \cite{Popov-short} and \cite{Zhu}. The case of sixteen dimensions was studied in \cite{Elashvili} and \cite{Harvey}. It is expected that beyond sixteen dimensions difficulties of principle arise. Another useful reference on the classification of spinors is the thesis \cite{Charlton}. 
Another, harder, version of the problem is the classification of the possible spinor orbits of the real spin group ${\rm Spin}(p,q)$. Significantly less is known about the real case, and much of the existing knowledge is summarised in \cite{Bryant}. 

In every dimension, there exists the orbit of smallest dimension. The corresponding spinors are known as \textbf{pure}, or in some literature simple. They have been studied since Cartan \cite{Cartan} and Chevalley \cite{Chevalley}. The other useful references are \cite{Benn-Tucker} and \cite{Budinich:1989bg}. There is a beautiful geometry related to pure spinors, and this will play an important role in what follows. 

Spinors are also representations of the Clifford algebra ${\rm Cl}(n)$, and every spinor has its associated null subspace, which is the subspace in the vector space that generates ${\rm Cl}(n)$ that annihilates the spinor. For pure spinors this subspace has the maximal possible dimension. For spinor that are not pure (which is the generic case in higher dimensions), the dimension of the annihilator subspace is not maximal, and in the generic case zero. A coarse classification of spinors according to their nullity has been given in \cite{TT}. Further, the classification of the possible types of pure spinors of the real spin group ${\rm Spin}(p,q)$ was given in \cite{KT}. 

In this paper we revisit the problem of classification of spinors. We develop a new approach that considers spinors as constructed from their pure spinor constituents. The case of odd dimensions is controlled by the next larger even dimension, and for this reason we restrict our attention to the case of even dimensions. In this paper we will only consider the problem over $\C$. Consider the Weyl irreducible (semi-spinor) representation $S^+$ of ${\rm Spin}(2d,\C)$. Our definitions of $S^\pm$, as well as the model for ${\rm Cl}(2d)$ that we use is described below. Given a pure spinor $\psi_0\in S^+$, let $M(\psi_0)$ be the corresponding annihilator subspace  
$$ \C^{2d} \supset M(\psi_0)=\{ u\in \C^{2d}: u\psi_0=0\},$$
where $u$ acts on $\psi_0$ by Clifford multiplication. Let us denote by $a_{1,\ldots,d}$ a basis in $M(\psi_0)$. Note that $M(\psi_0)$ is a totally null subspace $g(u,v)=0, \forall u,v\in M(\psi)$, as easily follows from the Clifford defining relation $uv+vu=2g(u,v)$. Here $g$ is the metric on $\C^{2d}$, which is the standard metric $g={\rm diag}(1,\ldots, 1)$ in any orthonormal basis. It is always possible (not in a unique way) to complete the vectors $a_{1,\ldots,d}$ by another set of $d$ null vectors that we denote by $a^\dagger_{1,\ldots,d}$ to form a basis of $\C^{2n}$ in which
\be\label{canonical-basis}
 g(a_i, a_j)=0, \quad g(a^\dagger_i, a^\dagger_j)=0, \quad g(a_i, a^\dagger_j)=\delta_{ij}, \qquad i,j=1,\ldots,d.
 \ee
We call such a basis of $C^{2d}$ a \textbf{canonical basis}.  
Given a canonical basis in $\C^{2d}$, a basis in $S^+$ is obtained by applying an \textbf{even number} of all possible (distinct) creation operators $a_i^\dagger$ to $\psi_0$
\be\label{spinor-basis}
S^+ ={\rm Span}(\psi_0,  \prod_{i<j} a_i^\dagger a_j^\dagger \psi_0 ,\prod_{i<j<k<l} a_i^\dagger a_j^\dagger a_k^\dagger a_l^\dagger\psi_0,\ldots \}.
\ee
It is easy to check that all $2^{d-1}$ vectors in the basis so constructed are pure spinors, and so every spinor can be represented as a linear combination of at most $2^{d-1}$ pure spinors. However, as is known since Cartan and Chevalley, the sum of two pure spinors can be a pure spinor, and so the description of a spinor as a sum of as many as $2^{d-1}$ pure spinors is superfluous. More precisely, let $\psi,\phi\in S^+$ be two pure spinors of the same helicity. Then the possible dimensions of the intersection $M(\psi)\cap M(\phi)$ are $d-2m, m=0,1,\ldots, [d/2]$. When $M(\psi),M(\phi)$ intersect in $d$ dimensions, the spinors $\psi,\phi$ are proportional to each other. When $M(\psi),M(\phi)$ intersect in $d-2$ dimensions, any linear combination of $\psi,\phi$ is still a pure spinor, see e.g. Proposition 5 in \cite{Budinich:1989bg}.

We approach the problem of classifying spinors that are linear combinations of the pure spinors from (\ref{spinor-basis}) combinatorially. Thus, let us consider a spinor $\psi$ that can be represented as a linear combination of $k$ pure spinors
\be\label{spinor-pure-sum}
\psi = \sum_{\alpha=1}^k c_\alpha \psi_\alpha . 
\ee
We get a tractable problem if we impose two restrictions: 
\begin{enumerate}
\item There exists a canonical (\ref{canonical-basis}) basis $a_i, a_i^\dagger$ of $C^{2d}$ such that every pure spinor $\psi_\alpha, \alpha=1,\ldots, k$ in (\ref{spinor-pure-sum}) is of the form of one of the basis vectors in (\ref{spinor-basis}); 
\item The pairwise intersections $M(\psi_\alpha)\cap M(\psi_\beta)$ of null subspaces of pure spinors that appear in (\ref{spinor-pure-sum}) are at most in $d-4$ dimensions, so that no sum of two pure spinors in (\ref{spinor-pure-sum}) is a pure spinor. 
\end{enumerate}
We represent the spinor $\psi$ in (\ref{spinor-pure-sum}) as a simplex with $k$ vertices labelled by $\psi_\alpha$, and we label the edges by the dimensions of the intersections $M(\psi_\alpha)\cap M(\psi_\beta)$. Our main conclusion is that, in dimension up to and including fourteen, there is only a limited number of distinct graphs that can be produced this way. We will call the graph corresponding to an impure spinor $\psi$ of the type (\ref{spinor-pure-sum}) {\it the combinatorial type of $\psi$.} Analysing the types of possible graphs quickly leads to the known classification results in dimensions up to and including fourteen. 

Another important point of our analysis is that, as we shall see, there is a new phenomenon in dimensions twelve that implies that not all distinct graphs correspond to distinct orbits. This phenomenon has already been described in \cite{Igusa}. Thus, it turns out that in dimension twelve there is an impure spinor of impurity four (i.e. one that is a linear combination of four pure spinors) that is in the same orbit as an impurity two spinor. This will be seen to imply that there are not only the edge intersection number constraints, but also the tetrahedral intersection number ones. Taking this phenomenon into account in higher dimensions allows to further restrict the number of possible combinatorial types of spinors. This allows us to reproduce all the known orbits in dimension up to and including fourteen by completely elementary combinatorial considerations. 

Our techniques are also applicable in dimension sixteen and higher. However, already in dimension sixteen our technology produces hundreds of different combinatorial types of spinors. Dimension sixteen is not treated here in details, see however some general remarks in Section \ref{sec:comb}, but we hope to return to this very interesting case in a separate publication. Dimension eighteen is even harder, but it is possible that by a further development of our method, and by utilising appropriate symbolic manipulation software, the classification of spinors using our method can also be achieved in this case. 

The organisation of this paper is as follows. We start by a brief description of the model for the Clifford algebra ${\rm Cl}(2d)$ that we use in Section \ref{sec:prelim}. We also explain here how the differential forms that arise as spinor bilinears are computed. The central section is \ref{sec:comb}, where we explain the combinatorics we use, and analyse for which dimensions this problem has a solution. The significance of dimension sixteen as one where new phenomena appear becomes clear here. We show here that there is bound on impurity even in dimensions as high as sixteen and eighteen, this bound coming from the tetrahedral intersection numbers considerations. To illustrate the power of the combinatorial approach we develop, we use it to reproduce the known classifications of spinors in dimensions up to and including twelve in \ref{sec:twelve} and in dimension fourteen in Section \ref{sec:fourteen}. 

\section{Creation/annihilation operators and a model for ${\rm Cl}(2d)$}
\label{sec:prelim}

\subsection{A model for the Clifford algebra}

A concrete model for ${\rm Cl}(2d)$ is obtained by selecting a pure spinor $\psi_0$. Let $a_i, i=1,\ldots, d$ be all null directions spanning the annihilator subspace $M(\psi_0)$. Let us complete them with a set of complementary null directions $a_i^\dagger$ which have the property $g(a_i, a_j^\dagger)=\delta_{ij}$. We will refer to $a_i$ as annihilation operators (because they annihilate $\psi_0$ when they act on it by the Clifford multiplication), and $a_i^\dagger$ as creation operators (because their action on $\psi_0$ produces other spinors). Consider the space
\be
S= {\rm Span}(\psi_0, a_i^\dagger \psi_0, \prod_{i<j} a_i^\dagger a_j^\dagger \psi_0, \ldots, a_1^\dagger \ldots a_d^\dagger \psi_0\},
\ee
spanned by all creation operators applied to $\psi_0$. It is clear that this space is isomorphic to the space of all differential forms in $\C^d$
\be
S \sim \Lambda^\bullet(\C^d).
\ee
Indeed, one just identifies
\be
 a_{i_1}^\dagger a_{i_2}^\dagger \ldots a_{i_k}^\dagger \psi_0 \sim e^{i_1 i_2 \ldots i_k}\in \Lambda^\bullet(\C^d).
\ee
Here $e^{i_1 i_2 \ldots i_k}\equiv e^{i_1}\wedge e^{i_2} \wedge \ldots \wedge e^{i_k}$, and $e^i, i=1,\ldots, d$ is a basis in $\Lambda^1(\C^d)$. 
Under this identification, spinors are identified with differential forms in $\Lambda^\bullet(\C^d)$.
The Clifford multiplication by $a_i$ is realised as the operator that removes the factor of $e^i$ from a spinor, and $a_i^\dagger$ is the operator that adds $e^i$. When written as an action on the basis vectors in $S$ this reads
\be
a_i^\dagger e^{i_1 i_2 \ldots i_k} = e^i\wedge e^{i_1 i_2 \ldots i_k}, \\ \nonumber
a_i e^{i_1 i_2 \ldots i_k} = \delta_i^{i_1}e^{i_2 \ldots i_k} - \delta_i^{i_2} e^{i_1 i_3 \ldots i_k}
+ \ldots + (-1)^{k-1} e^{i_1 i_2 \ldots i_{k-1}}.
\ee
 
\subsection{Lie algebra $\mathfrak{spin}(2d,\mathbb{C})$}

We can now give several useful descriptions of the Lie algebra $\mathfrak{spin}(2d,\mathbb{C})$. One description is to say that $\mathfrak{spin}(2d,\mathbb{C})$ is generated by all possible quadratic combinations of $a_i, a_i^\dagger$. A generic Lie algebra element can then be written as
\begin{equation}\label{lie algebra}
X=\frac{1}{2}A^{ij}(a_{i}a_{j}^{\dagger}- a_j^\dagger a_i) +\frac{1}{2}\beta^{ij}a_{i}^{\dagger}a_{j}^{\dagger}+\frac{1}{2} B^{ij}a_{i}a_{j}.
\end{equation}
Here the summation convention is implied, and $A^{ij}$ is an arbitrary $n\times n$ matrix, while $B^{ij}, \beta^{ij}$ are arbitrary anti-symmetric $n\times n$ matrices. 

Another description is to introduce linear combinations of $a_i, a_i^\dagger$ that have the property that the metric is diagonal on them
\be
\Gamma_{i}:=a_i+a_i^{\dagger}, \quad \Gamma_{i+d}:=\sqrt{-1}(a_i-a_i^{\dagger}). 
\ee
We then have $g(\Gamma_a,\Gamma_b)= 2\delta_{ab}, a,b=1,\ldots, 2d$. The Lie algebra $\mathfrak{spin}(2d,\mathbb{C})$ is then generated by all commutators 
\be
\mathfrak{spin}(2d,\mathbb{C}) = {\rm Span}( [\Gamma_{a},\Gamma_{b}]), \quad a,b=1,\ldots, 2d.
\ee

\subsection{Invariant inner product}

The formalism that identifies spinors with elements of $\Lambda^\bullet(\C^d)$ allows for a very simple explicit expression for the invariant inner product in $S$. Thus, the inner product between two spinors $\phi, \psi\in \Lambda^\bullet(\mathbb{C}^d)$ is given by
\begin{equation}\label{inner}
\langle\psi,\phi\rangle=\psi\wedge\sigma(\phi)\Big\vert_{\text{top}}.
\end{equation}
Here $\sigma(\phi)$ is an operation that acts on a basis $e^{i_1 i_2 \ldots i_k}\in \Lambda^\bullet(\mathbb{C}^d)$ by reordering it
\be
\sigma(e^{i_1 i_2 \ldots i_k}) = e^{i_k i_{k-1} \ldots i_1}.
\ee
The restriction to the top form selects the coefficient of $e^{12\ldots d}$. The fact that (\ref{inner}) is a $\mathfrak{spin}(2d,\mathbb{C})$-invariant inner product is checked by a computation. 

\subsection{Geometric map}

We now define a notion of a {\it geometric map}. This is a map that produces all possible differential forms from a spinor. The differential forms arise as spinor bilinears. Explicitly
\begin{equation}
    B_{k}(\psi,\psi):=\frac{1}{k!}\langle \psi, \Gamma_{b_1}\ldots\Gamma_{b_k} \psi\rangle E^{b_1}\wedge \ldots \wedge E^{b_k}.
\end{equation}
Here $E^{a}\in \Lambda^1(\C^{2d})$ is an orthonormal basis in the space of 1-forms in $\C^{2d}$. 

It is useful to recast this calculation as one in terms of creation/annihilation operators. We write
\begin{equation}
    E^i=\frac{e^i+\bar{e}^i}{2}, \ E^{i+d}=\frac{e^i-\bar{e}^i}{2\sqrt{-1}},
\end{equation}
where $e^i, \bar{e}^i \in \C^{2d}$ are null and $g(e_i,\bar{e}_j)=\delta_{ij}$. Then the geometric map takes the form
\begin{equation}
    B_k(\psi,\psi) =\frac{1}{k!} \langle \psi, (a_{i_1}e^{i_1}+a_{i_1}^{\dagger}\bar{e}^{i_1})\ldots(a_{i_k}e^{i_k}+a_{i_k}^{\dagger}\bar{e}^{i_k})\psi\rangle.
\end{equation}
This is the form of the geometric map that will be used for calculations below. To write this expression, we are making a distinction between creation/annihilation operators acting on spinors and 1-forms $e^i, \bar{e}^i\in \C^{2d}$. 

\section{Combinatorial problem}
\label{sec:comb}

Let us consider $k$ pure spinors $\psi_\alpha, \alpha=1,\ldots, k$ of the same helicity in $2d$ dimensions. We have assumed that there is a choice of basis $a_i, a_i^\dagger$ in $\C^{2d}$ that is canonical, see (\ref{canonical-basis}), and such that every pure spinor $\psi_\alpha$ is one of the basic pure spinors in (\ref{spinor-basis}). 

\subsection{Occupation numbers}

Let us think of every one of $2d$ basis vector $a_i, a_i^\dagger$ in $\C^{2d}$ as a box. Every one of the spinors $\psi_\alpha$ is a pure spinor whose null subspace $M(\psi_\alpha)$ is $d$-dimensional, and thus occupies precisely $d$ of the $2d$ available boxes. Let us introduce the box occupation numbers, and denote by 
\be
n_0, n_1, \ldots, n_k
\ee
the number of boxes occupied by $0, 1, \ldots, k$ pure spinors respectively. Each box is occupied by some number (possibly zero) of pure spinors, and so is counted. Therefore we have
\be\label{rel-1}
n_k + \ldots + n_0= 2d,
\ee
which is the total number of boxes. Another obvious relation is that each pure spinor occupies precisely $d$ of the boxes, and therefore
\be\label{rel-2}
k n_k + (k-1) n_{k-1} + \ldots n_1 = k d. 
\ee
Another relation between these occupation numbers is the duality. If a pure spinor occupies some box, it does not occupy the dual to it box. For example, if $a_i^\dagger$ is occupied then $a_i$ cannot be occupied by the same pure spinor, because the space spanned by the boxes occupied by some pure spinor must be totally null. This implies that there is a duality 
\be
n_k = n_0, \qquad n_{k-1} = n_1, \qquad {\rm etc.}
\ee
In view of this duality the two relations (\ref{rel-1}) and (\ref{rel-2}) are equivalent, and we only need one of them, e.g. (\ref{rel-1}). 
\begin{example}
Let us consider a simple example in dimension eight of a spinor that is not pure. This spinor is a sum of two pure spinors whose null subspaces do not intersect. Let one of these pure spinors, say $\psi_1$, be the vacuum state in (\ref{spinor-basis}) from which all other possible basis pure spinors are built. It is convenient to denote pure spinors by their corresponding null subspaces. And so we take
\be
\psi_1 = a_1 a_2 a_3 a_4.
\ee
This is somewhat of an abuse of notation, the precise meaning is that $M(\psi_1) = {\rm Span}(a_1, a_2, a_3, a_4)$. But the above is quicker to write and should not lead to any confusion. Let us take the other spinor to be 
\be
\psi_2= a_1^\dagger a_2^\dagger a_3^\dagger a_4^\dagger, 
\ee
where we again represent a pure spinor by its null subspace. In this case, $d=4$, and the occupation numbers are clearly $n_2=n_0=0$ and $n_1=8$. Indeed, there are no vectors (boxes) that are occupied by either zero or two pure spinors. And every box is occupied by precisely one pure spinor. Thus, all relations listed above are satisfied, which illustrates our combinatorial construction.
\end{example}
\begin{example} Let us consider a more involved example of a general spinor in eight dimensions given by a linear combination of all pure spinors appearing in (\ref{spinor-basis}). In this case $k=8$, as there are eight different possible states in (\ref{spinor-basis}). It is easy to see that every of the available 8 boxes is occupied precisely by 4 different pure spinors, so $n_4=8$. Again, all relations are satisfied. 
\end{example}

\subsection{Edge intersection numbers}

We now form a simplex with $k$ vertices and $k(k-1)/2$ edges. Each vertex represents a pure spinor, and each edge represents the intersection of the null subspaces of the two pure spinors that it connects. Let us introduce an index $I=1,\ldots, k(k-1)/2$ that labels the edges. Let us label the edge $I$ that connects vertices $\alpha,\beta$ with the dimension 
\be
e_I\equiv e_{\widehat{\alpha\beta}} = {\rm dim}( M(\psi_\alpha)\cap M(\psi_\beta))
\ee
of the intersection subspace of the null subspaces  $M(\psi_\alpha), M(\psi_\beta)$. As we know the possible dimensions are $d-2, d-4, \ldots$, where we assume that the pure spinors in the sum (\ref{spinor-pure-sum}) are distinct, and so their null subspaces cannot intersect in $d$ dimensions. But even the intersection number $d-2$ should not be allowed, as a pair of pure spinors intersecting in $d-2$ dimensions sums up to a pure spinor. Below we will see how to take this constraint into account in a useful way. 

There is a relation between the occupation numbers and the sum of the intersection numbers. We have
\be\label{occupation-intersection}
\frac{k(k-1)}{2} n_k + \frac{(k-1)(k-2)}{2} n_{k-1} + \ldots + n_2 = \sum_{I=1}^{k(k-1)/2} e_I.
\ee
Indeed, e.g. $n_k$ boxes that are occupied $k$ times each give rise to $k(k-1)/2$ pairs contributing to the sum of the intersection numbers. 

Let us illustrate the relation (\ref{occupation-intersection}) on the two examples already considered. 
\begin{example} For the spinor of impurity $k=2$ in eight dimensions the two pure constituents $\psi_1=a_1 a_2 a_3 a_4, \psi_2= a_1^\dagger a_2^\dagger a_3^\dagger a_4^\dagger$ have intersection number zero. So, we represent their linear combination of a simplex consisting of two vertices, and the edge connecting them. The edge has intersection number zero. The only non-vanishing occupation number is $n_1=8$, which does not appear in (\ref{occupation-intersection}). The right-hand side is also zero, as the only intersection number is zero in this case. We represent this configuration by the following diagram

\be
\begin{tikzpicture}
\coordinate [label=left:A] (A) at (-1,0);
\coordinate [label=right:B] (B) at (1,0);

\draw[blue, very thick] (A) -- (B);
\node[fill=white] (c) at ($(B)!0.5!(A)$) {$0$};

\end{tikzpicture} \qquad
\raise2ex\vbox{
\hbox{$\psi_A=a_1 a_2 a_3 a_4,$}
\hbox{$\psi_B=a_1^\dagger a_2^\dagger  a_3^\dagger a_4^\dagger .$}}
\ee

\end{example}
\begin{example} Let us now consider the spinor obtained as the general linear combination of all possible pure spinors. In this case we get a simplex with 8 vertices, and the occupation numbers are $n_4=8$. There are $28$ edges. From each vertex, there emanates an edge with intersection number zero, as well as six edges with intersection number two. This means that there are 24 edges with intersection number two. The relation (\ref{occupation-intersection}) is verified because we have $(4*3/2)*8=48$ on the left-hand side and $2*24$ on the right-hand side. We refrain from drawing this graph, as it is discussed here only for illustration purposes, and is not needed below. 
\end{example}

\subsection{Tetrahedral intersection numbers}

We have labelled each edge by the dimension of the null subspace that is common to both vertices bounded by this edge. This can be continued, and one can introduce similar intersection numbers for each triangle, tetrahedron, etc. As will become clear below, the tetrahedral intersection number will play an important role. 

There are $k(k-1)(k-2)(k-3)/24$ tetrahedra in the complete graph built on $k$ vertices. We label each of these tetrahedra by the dimension of the common null subspace 
\be
t_{\widehat{\alpha\beta\gamma\delta}} = {\rm dim}( M(\psi_\alpha)\cap M(\psi_\beta)\cap M(\psi_\gamma)\cap M(\psi_\delta)). 
\ee
We now have a relation similar to (\ref{occupation-intersection})
\be\label{occupation-intersection-4}
\frac{k(k-1)(k-2)(k-3)}{24} n_k + \frac{(k-1)(k-2)(k-3)(k-4)}{24} n_{k-1} + \ldots + n_4 = \sum_{\widehat{\alpha\beta\gamma\delta}} t_{\widehat{\alpha\beta\gamma\delta}}.
\ee
The reason for why this relation is useful is that the tetrahedral intersection numbers cannot be too high. Thus, as we shall see below, we must have $t\leq d-7$ in order to obtain a linear combination of pure spinors that cannot be further reduced. 

\begin{example} We will illustrate this formula by the spinor that is given by the linear combination of all basic pure spinors in six dimensions. The number of basis vectors in (\ref{spinor-basis}) is 4. All the null directions are shared by exactly two pure spinors, so $n_2=6$. The intersection dimension on each edge is $e=1$. There is no directions shared by all four pure spinors, and so the tetrahedral intersection number for the single tetrahedron is $t=0$.
\end{example}
\begin{example} Another example is the complete set of pure spinors in eight dimensions. There are eight pure spinors and $n_4=8$. There are in total $70$ tetrahedra, but most have intersection number zero. There are eight tetrahedra with intersection number one, and they are given by
\be\nonumber
a_1 a_2 a_3 a_4, \quad a_2^\dagger a_3^\dagger a_1 a_4, \quad a_2^\dagger a_4^\dagger a_1 a_3, \quad a_3^\dagger a_4^\dagger a_1 a_2, \quad \text{all share direction $a_1$} \\ \nonumber
a_1 a_2 a_3 a_4, \quad a_1^\dagger a_3^\dagger a_2 a_4, \quad a_1^\dagger a_4^\dagger a_2 a_3, \quad a_3^\dagger a_4^\dagger a_1 a_2, \quad \text{all share direction $a_2$} \\ \nonumber
a_1 a_2 a_3 a_4, \quad a_1^\dagger a_2^\dagger a_3 a_4, \quad a_1^\dagger a_4^\dagger a_2 a_3, \quad a_2^\dagger a_4^\dagger a_1 a_3, \quad \text{all share direction $a_3$} \\ \nonumber
a_1 a_2 a_3 a_4, \quad a_1^\dagger a_2^\dagger a_3 a_4, \quad a_1^\dagger a_3^\dagger a_2 a_4, \quad a_2^\dagger a_3^\dagger a_1 a_4, \quad \text{all share direction $a_4$}  \\ \nonumber
a_1^\dagger a_2^\dagger a_3^\dagger a_4^\dagger, \quad a_1^\dagger a_2^\dagger a_3 a_4, \quad a_1^\dagger a_3^\dagger a_2 a_4, \quad a_1^\dagger a_4^\dagger a_2 a_3, \quad \text{all share direction $a_1^\dagger$}
\\ \nonumber
a_1^\dagger a_2^\dagger a_3^\dagger a_4^\dagger, \quad a_1^\dagger a_2^\dagger a_3 a_4, \quad a_2^\dagger a_3^\dagger a_1 a_4, \quad a_2^\dagger a_4^\dagger a_1 a_3, \quad \text{all share direction $a_2^\dagger$}
\\ \nonumber
a_1^\dagger a_2^\dagger a_3^\dagger a_4^\dagger, \quad a_1^\dagger a_3^\dagger a_3 a_4, \quad a_2^\dagger a_3^\dagger a_1 a_4, \quad a_3^\dagger a_4^\dagger a_1 a_2, \quad \text{all share direction $a_3^\dagger$}
\\ \nonumber
a_1^\dagger a_2^\dagger a_3^\dagger a_4^\dagger, \quad a_1^\dagger a_4^\dagger a_2 a_3, \quad a_2^\dagger a_4^\dagger a_1 a_3, \quad a_3^\dagger a_4^\dagger a_1 a_2, \quad \text{all share direction $a_4^\dagger$}
\ee
Thus, there are eight tetrahedra with $t=1$, and (\ref{occupation-intersection-4}) holds. 
\end{example}

We will now analyse what all the equations imply together with the constraints on the edge and tetrahedral intersection numbers. We consider the case of $k$ odd and even separately.

\subsection{Odd number of pure spinors}

When $k+1$ is even all occupation numbers are in pairs related by the duality. We have the single equation
\be\label{odd}
n_k + n_{k-1} + \ldots + n_{(k+1)/2} = d.
\ee
We also rewrite the relation (\ref{occupation-intersection}) in terms of the higher occupation numbers
\be
\frac{k(k-1)}{2} n_k + \frac{(k-1)(k-2)}{2} n_{k-1} +  \left( \frac{(k-2)(k-3)}{2} +1\right) n_{k-2} + \ldots
\\ \nonumber
+ \left( \frac{(k-a)(k-a-1)}{2} + \frac{a(a-1)}{2}\right) n_{k-a} + \ldots +
\frac{(k-1)^2}{4} n_{(k+1)/2} = \sum_{I=1}^{k(k-1)/2} e_I.
\ee
Here $a$ goes up to the maximal value $a_{max}=(k-1)/2$. 
If we substitute
\be\label{lowest-odd}
n_{(k+1)/2} = d - n_k - n_{k-1} - \ldots -n_{(k+3)/2}.
\ee
we get
\be
\left( \frac{k(k-1)}{2} -\frac{(k-1)^2}{4} \right) n_k + \left( \frac{(k-1)(k-2)}{2} - \frac{(k-1)^2}{4}\right) n_{k-1} +\ldots
\\ \nonumber 
+  \left( \frac{(k-a)(k-a-1)}{2} +\frac{a(a-1)}{2} - \frac{(k-1)^2}{4} \right) n_{k-a} + \ldots
= \sum_{I=1}^{k(k-1)/2} e_I - \frac{(k-1)^2}{4}  d.
\ee
The $n_{(k+1)/2}$ term has been subtracted, so the last term on the left-hand side is one containing $n_{(k+3)/2}$. One can check that the coefficients in front of all the terms on the left-hand side are positive, and so left-hand side is a non-negative number. On the other hand, as we have already discussed, the maximal value of the intersection number in $2d$ dimensions is $d-4$. This is so that no sum of the pure spinor is pure. And so the right-hand side of the above is smaller or equal than
\be
\frac{k(k-1)}{2} (d-4) - \frac{(k-1)^2}{4} d =\frac{(k-1)}{4}( (k+1)d-8k).
\ee
This must be non-negative, because the sum on the left-hand side is non-negative, which gives, for $k>1$
\be\label{ineq-odd}
k(d-8) +d \geq 0.
\ee
It is clear that for $d=2, 3, 4, 5$ the only possible solution (recall that here $k$ is odd by assumption) is $k=1$. For $d=6$ we get another possible solution $k=3$, while for $d=7$ the solutions are $k=1,3,5,7$. For $d\geq 8$ the inequality is always satisfied, and no restriction on the number of pure spinors come this way. Thus, our combinatorial analysis confirms that spinors in dimension sixteen behave very differently to the situation in lower dimensions. 

Let us now analyse the consequences of the restriction on the tetrahedral intersection numbers. We use the duality and rewrite (\ref{occupation-intersection-4}) as
\be\nonumber
\frac{k(k-1)(k-2)(k-3)}{24} n_k + \frac{(k-1)(k-2)(k-3)(k-4)}{24} n_{k-1} + \frac{(k-2)(k-3)(k-4)(k-5)}{24} n_{k-2}
\\ \nonumber
+ \frac{(k-3)(k-4)(k-5)(k-6)}{24} n_{k-3}+\left( \frac{(k-4)(k-5)(k-6)(k-7)}{24} +1\right) n_{k-4} + \ldots
\\ \nonumber
+ \left( \frac{(k-a)(k-a-1)(k-a-2)(k-a-3)}{24} +\frac{a(a-1)(a-2)(a-3)}{24} \right) n_{k-a} + \ldots
\\ \nonumber
+ \frac{(k-1)(k-3)^2(k-5)}{96} n_{(k+1)/2} = \sum_{\widehat{\alpha\beta\gamma\delta}} t_{\widehat{\alpha\beta\gamma\delta}}.
\ee
We then again use (\ref{lowest-odd}) and take leave all terms containing the occupation numbers on the left-hand side, while taking the term containing $d$ to the right. We get the following right-hand side
\be
\sum_{\widehat{\alpha\beta\gamma\delta}} t_{\widehat{\alpha\beta\gamma\delta}} - \frac{(k-1)(k-3)^2(k-5)}{96} d.
\ee
This must be non-negative. As will become clear from our analysis below, the largest value of each tetrahedral intersection number that does not lead to a reduction in the number of pure spinors is $t\leq d-7$. The number of tetrahedra is $k(k-1)(k-2)(k-3)/24$. So, the following inequality must be satisfied
\be
\frac{k(k-1)(k-2)(k-3)}{24} (d-7)  - \frac{(k-1)(k-3)^2(k-5)}{96} d \geq 0.
\ee
For $k>3$ this is equivalent to
\be
4k(k-2)(d-7) - (k-3)(k-5) d\geq 0,
\ee
or
\be
k^2(28-3d) - 56 k + 15d \leq 0.
\ee
This does give interesting restrictions. We don't need to consider the dimensions up to and including $d=6$, because we already know that the maximal value of $k$ in $d=6$ is $k=3$. The new constraints arise in $d=7$, where we get that $k\leq 5$. This is a constraint stronger than the one we had previously by considering only the edge intersection numbers. The new case is $d=8$, where we had no constraint coming from the edge intersection numbers. In this case we have an equivalent inequality $4k^2-56k +120\leq 0$. This must be satisfied by odd integers $k>3$. The largest such integer is $k=11$, which gives us a bound on the number of pure spinors in the case of sixteen dimensions. This bound arises solely from considering the restriction on the tetrahedral intersection numbers. 

\subsection{Even number of pure spinors}

We now perform a similar analysis for the even number of pure spinors. We have 
\be\label{even}
2n_k + 2n_{k-1} + \ldots + 2 n_{k/2+1}+ n_{k/2} = 2d.
\ee
The expression for the sum of the intersection numbers in terms of the independent occupation numbers is now
\be
\frac{k(k-1)}{2} n_k + \frac{(k-1)(k-2)}{2} n_{k-1} +  \left( \frac{(k-2)(k-3)}{2} +1\right) n_{k-2} + \ldots
\\ \nonumber
+ \left( \frac{(k-a)(k-a-1)}{2} + \frac{a(a-1)}{2}\right) n_{k-a} + \ldots 
+ \frac{k(k-2)}{8} n_{k/2} = \sum_{i=1}^{k(k-1)/2} e_i.
\ee
Expressing $n_{k/2}$ in terms of the other occupation numbers, taking the factor of $d k(k-2)/ 4$ to the left-hand side, and using the fact that the maximal possible value of the edge intersection numbers is $e\leq (d-4)$, we get the maximal possible value of the right-hand side to be
\be
\frac{k(k-1)}{2} (d-4) - \frac{k(k-2)}{4} d= \frac{k}{4}(k (d-8) + 8).
\ee
Thus, we get instead
\be
k(d-8) + 8 \geq 0.
\ee
For $k\leq 7$ this implies
\be
k \leq \frac{8}{8-d}.
\ee
For $d=4,5$ this gives $k\leq 2$, which tells us that there the only impure spinors are those consisting of two pure spinors in these numbers of dimensions. This is a known result, which we have reproduced by our method. For $d=6$ we get $k\leq 4$. For $d=7$ we get $k\leq 8$. For $d\geq 8$ there are no constraints arising on $k$ this way solely from the edge intersection number considerations. 

Let us now study the consequences of the restriction on the maximal possible value of the tetrahedral intersection number. We first rewrite (\ref{occupation-intersection-4}) using the duality
\be\nonumber
\frac{k(k-1)(k-2)(k-3)}{24} n_k + \frac{(k-1)(k-2)(k-3)(k-4)}{24} n_{k-1} + \frac{(k-2)(k-3)(k-4)(k-5)}{24} n_{k-2}
\\ \nonumber
+ \frac{(k-3)(k-4)(k-5)(k-6)}{24} n_{k-3}+\left( \frac{(k-4)(k-5)(k-6)(k-7)}{24} +1\right) n_{k-4} + \ldots
\\ \nonumber
+ \left( \frac{(k-a)(k-a-1)(k-a-2)(k-a-3)}{24} +\frac{a(a-1)(a-2)(a-3)}{24} \right) n_{k-a} + \ldots
\\ \nonumber
+ \frac{k (k-2)(k-4)(k-6)}{192} n_{k/2} = \sum_{\widehat{\alpha\beta\gamma\delta}} t_{\widehat{\alpha\beta\gamma\delta}}.
\ee
We then substitute the value of $n_{k/2}$ from (\ref{even}), and take the term containing the dimension $2d$ to the right-hand side. We use the fact that the maximal value of $t$ for each tetrahedron is $t=d-7$. This gives the following inequality that must be satisfied 
\be
\frac{k(k-1)(k-2)(k-3)}{24} (d-7)  - \frac{k(k-2)(k-4)(k-6)}{96} d \geq 0.
\ee
For $k>2$ this is equivalent to
\be
4(k-1)(k-3) (d-7) - d(k-4)(k-6) \geq 0,
\ee
or equivalently
\be
k^2(28-3d) - k(112-6d) + 12d + 84 \leq 0. 
\ee
For $d=7$ this gives $k\leq 6$. For $d=8$ this gives $k\leq 12$. Both give interesting and useful bounds on the numbers of possible pure spinors making up an impure spinor. 

\section{Dimension up to and including twelve}
\label{sec:twelve}

We now explicitly solve the combinatorial problem in each given dimension, listing the possible graphs allowed by the combinatorial analysis, finding a possible spinor representing each graph, and finding the (simple parts of the) stabiliser. In this section we reproduce the known classification of orbits up to and including twelve dimensions. 

\subsection{Action of a Cartan subgroup}

Our spinors are arbitrary linear combinations (\ref{spinor-pure-sum}) of pure spinors. It is useful to be able determine whether the coefficients in this linear combination can be scaled away by the action of the group. To this end, we consider the action of a suitably chosen Cartan subgroup. A convenient choice for the Cartan subalgebra is
\be
 \mathfrak{h}={\rm Span}\{t_i, i=1,\ldots,d\}, \qquad t_i := a_ia_i^{\dagger}-a_i^{\dagger}a_i.
 \ee
 Each of the pure spinors in the basis (\ref{spinor-pure-sum}) is of the form
 \be
 a_{i_1}^\dagger \ldots a_{i_k}^\dagger \psi_0.
 \ee
 Each such pure spinor is an eigenstate of each $t_i$. The eigenvalues are $+1$ if $i$ is not among $i_1,\ldots, i_k$ and $-1$ otherwise. This means that the eigenvalues of $e^{\lambda_i t_i}$ (no summation implied) are $e^{\lambda_i}$ if $i$ is not among $i_1,\ldots, i_k$ and $e^{-\lambda_i}$ otherwise. For every impure spinor we consider below, we will attempt to use the action of the Cartan generators to rescale all the coefficients $c_\alpha$ to make them equal. Whether this is possible or not of course depends on an impure spinor under consideration. 
  
 \subsection{Eight dimensions }
 
It is well-known that every Weyl spinor in dimensions two, four, and six is pure, with stabiliser ${\rm SU}(d)$. The first non-trivial problem arises in dimension eight. From the general analysis above, we expect to be able to construct an impure spinor of impurity two in this case. We are thus looking at $k=2$. We have $2n_2+ n_1= 8$, and also $n_2=e$, where $e$ is the intersection number between the two pure spinors. But given that the largest value of $e$ is $d-4=0$, we must have $n_2=0$. So, the only solution is given by two pure spinors whose null subspaces are complimentary. This is the configuration already considered as an example in the previous section
\be
\begin{tikzpicture}
\coordinate [label=left:A] (A) at (-1,0);
\coordinate [label=right:B] (B) at (1,0);

\draw[blue, very thick] (A) -- (B);
\node[fill=white] (c) at ($(B)!0.5!(A)$) {$0$};

\end{tikzpicture} \qquad
\raise2ex\vbox{
\hbox{$\psi_A=a_1 a_2 a_3 a_4,$}
\hbox{$\psi_B=a_1^\dagger a_2^\dagger  a_3^\dagger a_4^\dagger .$}}
\ee

The most general impure spinor arising in eight dimensions is thus $\psi=c_A \psi_A + c_B \psi_B$. However, it is clear that the coefficients $c_A, c_B$ can be made equal by the action of the Cartan subgroup. Indeed, we can simply take all $\lambda_i, i=1,\ldots,4$ to be equal in this case, and act with $t_1 t_2 t_3 t_4$. The action on $\psi_A$ gives $\psi_A \to e^{4\lambda}\psi_A$, while the action on $\psi_B$ gives $\psi_B \to e^{-4\lambda} \psi_B$. It is clear that there is a choice of $\lambda$ that makes $c_A=c_B$. Thus, there is just an overall constant factor in this spinor, that needs no longer to be considered. So, it is sufficient to analyse $\psi=\psi_A+\psi_B$.

To understand the stabiliser of an impure spinor, it is sufficient to build all possible spinor bilinears $B_k(\psi)$. The stabiliser of $\psi$ is then the subgroup of ${\rm SO}(2d)$ that stabilises all the differential forms $B_k(\psi)$. For a pure spinor in eight dimensions only $B_4(\psi_{A,B})\not=0$. Each of these 4-forms is decomposable and is given by the product of the null directions spanning the null subspace of $\psi_{A,B}$
\be
B_4(\psi_A) = \bar{e}_1 \wedge \bar{e}_2 \wedge \bar{e}_3 \wedge \bar{e}_4, \quad
B_4(\psi_B) = e_1 \wedge e_2 \wedge e_3 \wedge e_4.
\ee
Here we used the symbol of the wedge product to indicate that the directions $a_i, a_i^\dagger$ are treated as objects in $\Lambda^1(\C^{2d})$.
Given two pure spinors as above, we also have the following non-vanishing spinor bilinears 
\be
B_0(\psi_A, \psi_B) = 1, \quad B_2(\psi_A, \psi_B) = \omega, \quad B_4(\psi_A, \psi_B) =\frac{1}{2} \omega\wedge \omega,
\ee
where
\be
\omega = \bar{e}_1 \wedge e_1 + \bar{e}_2 \wedge e_2 + \bar{e}_3 \wedge e_3 + \bar{e}_4 \wedge e_4.
\ee
This shows that 
\be
B_4(\psi) =   \bar{e}_1 \wedge \bar{e}_2 \wedge \bar{e}_3 \wedge \bar{e}_4 +  e_1 \wedge e_2 \wedge e_3 \wedge e_4 + \frac{1}{2}  \omega\wedge \omega.
\ee
This is a 4-form in 8D known as the Cayley form. Its stabiliser is ${\rm Spin}(7)$, and so this is also the stabiliser of $\psi$. 

\subsection{Ten dimensions}

Again, we know that the only possible impure spinor is one with $k=2$. We have $2n_2+ n_1= 10$ and $n_2=e$. The only possible value for $e$ are $d-4=1$ (the intersection number must take odd values for odd $d$). So, the only possible solution is $n_2=1, n_1=8$. A possible solution representing this is the following configuration
\be
\begin{tikzpicture}
\coordinate [label=left:A] (A) at (-1,0);
\coordinate [label=right:B] (B) at (1,0);

\draw[blue, very thick] (A) -- (B);
\node[fill=white] (c) at ($(B)!0.5!(A)$) {$1$};

\end{tikzpicture} \qquad
\raise2ex\vbox{
\hbox{$\psi_A=a_1 a_2 a_3 a_4 a_5,$}
\hbox{$\psi_B=a_1^\dagger a_2^\dagger  a_3^\dagger a_4^\dagger a_5.$}}
\ee
We have chosen the null direction common to both pure spinors to be $a_5$. It is clear that the spinor given by a linear combination of these two pure spinors is effectively the unique impure spinor in 8 dimensions that was already described above. The additional dimension does not play any role. The rescaling of this spinor to one with an overall scale factor by the action of the Cartan subgroup is identical to the already considered 8D case. Below we will see that there is a similar impure spinor of impurity two in any dimension. Its stabiliser contains ${\rm Spin}(7)$ stabilising the Cayley form in 8D.

\subsection{12 Dimensions}

We expect that in this case $k\leq 3$. But we will also need to consider the case $k=4$, as the elimination of this case as not independent is what leads to the bound on the tetrahedral intersection number that played an important role in the previous section. 

\subsubsection{$k=2$}

The equation relating the occupation numbers is $2n_2+n_1 = 12$. We also have $n_2=e$, and the possible dimensions of the intersections are $d-4,d-6$ and thus $2,0$. When the intersection dimension is $e=2$ we have $n_2=n_0=2, n_1=8$. A possible representative is the configuration
\be
\begin{tikzpicture}
\coordinate [label=left:A] (A) at (-1,0);
\coordinate [label=right:B] (B) at (1,0);

\draw[blue, very thick] (A) -- (B);
\node[fill=white] (c) at ($(B)!0.5!(A)$) {$2$};

\end{tikzpicture} \qquad
\raise2ex\vbox{
\hbox{$\psi_A=a_1 a_2 a_3 a_4 a_5 a_6,$}
\hbox{$\psi_B=a_1^\dagger a_2^\dagger  a_3^\dagger a_4^\dagger a_5 a_6.$}}
\ee
This is again effectively the impure spinor in 8D. The simple part of the stabiliser is ${\rm Spin}(7)\times{\rm SL}(2)$. The spinor $\psi= c_A \psi_A + c_B \psi_B$ can again be rescaled by the action of the Cartan subgroup into the form of on overall coefficient times $\psi= \psi_A+\psi_B$.

To understand the stabiliser of this impure spinor, we compute the differential forms arising as the result of the geometric map. The differential forms $B_k(\psi,\psi)$ in 12D are only possible with $k=2, 6$. We get $B_2(\psi,\psi)= \bar{e}_5\wedge \bar{e}_6$, and $B_6(\psi,\psi)$ equals the Cayle form in the 8D space spanned by $e_{1,2,3,4}, \bar{e}_{1,2,3,4}$, times $\bar{e}_5\wedge \bar{e}_6$. This shows that the simple part of the stabiliser is indeed ${\rm Spin}(7)$ acting in the 8D subspace, and then ${\rm SL}(2)$ mixing the null directions $\bar{e}_5,\bar{e}_6$.

When the dimension of the intersection is zero, we have $n_1=12$ and the pure spinors are complementary
\be\label{12d-p2-2}
\begin{tikzpicture}
\coordinate [label=left:A] (A) at (-1,0);
\coordinate [label=right:B] (B) at (1,0);

\draw[blue, very thick] (A) -- (B);
\node[fill=white] (c) at ($(B)!0.5!(A)$) {$0$};

\end{tikzpicture} \qquad
\raise2ex\vbox{
\hbox{$\psi_A=a_1 a_2 a_3 a_4 a_5 a_6,$}
\hbox{$\psi_B=a_1^\dagger a_2^\dagger  a_3^\dagger a_4^\dagger a_5^\dagger a_6^\dagger.$}}
\ee
Considering the action of the Cartan subgroup with all $\lambda_i$ equal, it is again easy to see that the coefficients in the linear combination $c_A\psi_A + c_B \psi_B$ can be made equal, and so it is sufficient to consider $\psi=\psi_A+\psi_B$. 

To understand the stabiliser, we compute
\be
B_2(\psi)= \bar{e}_1 \wedge e_1 + \ldots \bar{e}_6 \wedge e_6.
\ee
The 6-form contains decomposable pieces comprising the directions $e_{1,\ldots 6}$ and $\bar{e}_{1,\ldots,6}$, as well as the cube of the 2-form $B_2(\psi)$. The stabiliser in this case contains ${\rm SL}(6)$ that mixes the null directions $e_i \to g_{ij} e_j$ and $\bar{e}_i \to g^{-1}_{ji} \bar{e}_j$, so that $B_2(\psi)$ remains invariant. 

\subsubsection{$k=3$}

Now consider possible solutions arising with a triple of pure spinors. The occupation numbers are related via $n_3+ 
n_2 =6$. If we denote the intersection dimensions by $e_1, e_2, e_3$, we have
\be\label{intersections3-12}
3 n_3 + n_2 = e_1+e_2+e_3.
\ee
Using the relation between $n_2,n_3$ this can be rewritten as
\be
n_3 = \frac{1}{2}\left(\sum_{I=1}^3 e_I-6\right).
\ee 
At the same time, the maximal possible dimension of the intersection in this number of dimensions is $e=2$. This means that there is a unique solution in the case of $k=3$, which is with all intersection numbers on edges having value $e_{1,2,3}=2$. This gives $n_3=n_0=0$, and then $n_2=n_1=6$. A possible representative is given by the following configuration
\be\label{12d-p3}
\begin{tikzpicture}
\coordinate [label=left:A] (A) at (-2*0.75,0);
\coordinate [label=above:B] (B) at (0,3.46*0.75);
\coordinate [label=right:C] (C) at (2*0.75,0);

\draw[blue, very thick] (A) -- (B);
\node[fill=white] (c) at ($(B)!0.5!(A)$) {$2$};

\draw[blue, very thick] (B) -- (C);
\node[fill=white] (c) at ($(C)!0.5!(B)$) {$2$};

\draw[blue, very thick] (A) -- (C);
\node[fill=white] (c) at ($(C)!0.5!(A)$) {$2$};
\end{tikzpicture} \qquad
\raise6ex\vbox{
\hbox{$\psi_A=a_1 a_2 a_3 a_4 a_5 a_6,$}
\hbox{$\psi_B=a_1^\dagger a_2^\dagger a_3^\dagger a_4^\dagger a_5 a_6,$}
\hbox{$\psi_C=a_1 a_2  a_3^\dagger a_4^\dagger a_5^\dagger a_6^\dagger.$}}
\ee

To show that the individual coefficients $c_{A,B,C}$ can all be made equal by the action of the Cartan subgroup, we consider the action of $e^{\lambda_1 t_1}$ and $e^{\lambda_5 t_5}$. We then want
\be
c=c_A e^{\lambda_1+\lambda_5}= c_B e^{\lambda_5-\lambda_1} = c_C e^{\lambda_1 - \lambda_5},
\ee
which is clearly possible. So, it is sufficient to consider $\psi=\psi_A+\psi_B+\psi_C$. 

To understand the stabiliser of this spinor, we compute $B_2(\psi)$. It is given by 3 terms, each being the wedge product of two of the null directions shared by each pair of pure spinors
\be
B_{2}(\psi)=\bar{e}_5\wedge \bar{e}_6+ \bar{e}_1 \wedge \bar{e}_2+ e_3 \wedge e_4. 
\ee
It is clear that the simple part of the stabiliser of $\psi$ is given by ${\rm Sp}(6)$, which is the symplectic group in the 6-dimensional (totally null) space spanned by $\bar{e}_{1,2}, e_{3,4}, \bar{e}_{5,6}$. 

\subsubsection{$k=4$}

The relation between the occupation numbers is $2n_4+ 2n_3+ n_2=12$. The relation between the occupation numbers and the intersection dimensions on the edges of the resulting tetrahedron is
\be
6 n_4 + 3 n_3 + n_2 = \sum_{I=1}^6 e_I.
\ee
We now express $n_2$ via $n_4, n_3$ and substitute into this equation to get
\be
4n_4+n_3= \sum_{i=1}^6 e_i -12,
\ee
The maximal possible value each intersection dimension is $e=2$, and so the maximal value of the sum of the intersection numbers is 12. We thus see that there is a single solution of the combinatorial problem in this case
\be
n_4=n_0=0, \qquad n_3=n_1 =0, \qquad n_2=12.
\ee
A possible representative of this solution is the following set of four pure spinors
\be\label{12d-p4}
\begin{tikzpicture}
\coordinate [label=left:A] (A) at (-1.5,0);
\coordinate [label=above:B] (B) at (0,1.5);
\coordinate [label=right:C] (C) at (1.5,0);
\coordinate [label=below:D] (D) at (0,-1.5);

\draw[blue, very thick] (A) -- (B);
\node[fill=white] (c) at ($(B)!0.5!(A)$) {$2$};

\draw[blue, very thick] (B) -- (C);
\node[fill=white] (c) at ($(C)!0.5!(B)$) {$2$};

\draw[blue, very thick] (C) -- (D);
\node[fill=white] (c) at ($(D)!0.5!(C)$) {$2$};

\draw[blue, very thick] (D) -- (A);
\node[fill=white] (c) at ($(A)!0.5!(D)$) {$2$};

\draw[blue, very thick] (A) -- (C);
\node[fill=white] (c) at ($(C)!0.75!(A)$) {$2$};

\draw[blue, very thick] (B) -- (D);
\node[fill=white] (c) at ($(D)!0.75!(B)$) {$2$};
\end{tikzpicture}\qquad
\raise6ex\vbox{
\hbox{$\psi_A=a_1 a_2 a_3 a_4 a_5 a_6,$ }
\hbox{$\psi_B=a_1^\dagger a_2^\dagger a_3^\dagger a_4^\dagger a_5 a_6,$}
\hbox{$\psi_C=a_1 a_2 a_3^\dagger a_4^\dagger a_5^\dagger a_6^\dagger,$}
\hbox{$\psi_D=a_1^\dagger a_2^\dagger a_3 a_4 a_5^\dagger a_6^\dagger.$}}
\ee

To understand the geometry of this orbit, we compute $B_2(\psi)$, where
\be
\psi= \sum_{\alpha=1}^4 c_\alpha \psi_\alpha.
\ee
We could have scaled away the individual coefficients using the action of the Cartan subgroup, but given that we want to reduce this case to the previously considered one, this may as well be done with general coefficients. Given that $\psi$ is a linear combination of 4 pure spinors, and $B_2$ vanishes for each pure spinor (only $B_6$ is different from zero), the non-zero contributions to $B_2(\psi)$ come from pairs of different pure spinors. However, each such pair intersects precisely in two dimensions, and $B_2(\psi_\alpha,\psi_\beta)$, for two pure spinors that intersection in two dimensions, is a multiple of the wedge product of the corresponding null directions. It is then clear that $B_2(\psi)$ is given by  the following six terms
\be
B_2(\psi)=c_A c_C \bar{e}_1 \wedge \bar{e}_2 + c_A c_D \bar{e}_3 \wedge \bar{e}_4 + c_A c_B \bar{e}_5 \wedge \bar{e}_6 \\ \nonumber
+ c_B c_D e_1 \wedge e_2 +c_B c_C e_3\wedge e_4 +c_C c_D e_5\wedge e_6.
\ee
Let us consider the terms involving only the directions $1,2$, i.e., $c_A c_C \bar{e}_1\wedge \bar{e}_2 + c_B c_D e_1\wedge e_2$. We want to show that there is a different canonically normalised null basis $b_1, b_2, \bar{b}_1, \bar{b}_2$ in the space spanned by $\bar{e}_{1,2}, e_{1,2}$, such that
\be\label{desired}
c_A c_C \bar{e}_1\wedge \bar{e}_2 + c_B c_D e_1\wedge e_2=\lambda (\bar{b}_1 \wedge b_1+ \bar{b}_2 \wedge b_2),
\ee
where $\lambda$ is some constant. We define
\be
\bar{b}_1 = \alpha \bar{e}_1 +\beta e_2, \quad b_1 = \gamma \bar{e}_2 + \delta e_1, \\ \nonumber
\bar{b}_2= \alpha \bar{e}_2 - \beta e_1, \quad b_2 = \delta e_2 - \gamma \bar{e}_1,
\ee
which satisfies 
\be
g(\bar{b}_1, \bar{b}_2)=g(b_1, b_2)=0, \qquad g(\bar{b}_1, b_1) = g(\bar{b}_2, b_2)=\alpha\delta +\beta \gamma.
\ee
It is then easy to check that the 2-form $\bar{b}_1 \wedge b_1 + \bar{b}_2 \wedge b_2$, does not have any $\bar{e}_1 \wedge e_1, \bar{e}_2\wedge  e_2$ terms when $\alpha\delta=\beta\gamma$. When this is satisfied
\be
\bar{b}_1 \wedge b_1 + \bar{b}_2 \wedge b_2= 2 \alpha\gamma \bar{e}_1\wedge \bar{e}_2 - 2\beta\delta e_1\wedge e_2.
\ee
We then choose $\delta = (2\alpha)^{-1}, \gamma=(2\beta)^{-1}$ to have the canonical normalisations. Then (\ref{desired}) holds for 
\be
\frac{\alpha^2}{\beta^2} = - \frac{c_A c_C}{c_B c_D}, \qquad \lambda^2 = - c_A c_B c_C c_D. 
\ee
Thus, (\ref{desired}) is indeed possible for an appropriate choice of a canonically normalised null basis in the space spanned by $\bar{e}_1, \bar{e}_2, e_1, e_2$. This means that there exists another canonical basis $\bar{b}_{1,\ldots,6}, b_{1,\ldots,6}$ in $\C^{12}$ such 
\be
B_2(\psi) \sim \bar{b}_1 \wedge b_1 + \ldots + \bar{b}_6 \wedge b_6 .
\ee
This shows that the impurity four spinor in twelve dimensions is in the same ${\rm Spin}(12)$ orbit as the impurity two spinor (\ref{12d-p2-2}) with the stabiliser ${\rm SL}(6)$. This 12D phenomenon is then at the root of some additional reductions that occur in dimensions 14D and higher, and reduce the number of possible combinatorial types of orbits. 

\section{14 Dimensions}
\label{sec:fourteen}

We now consider the more difficult case of fourteen dimensions. As we know, just taking the edge intersection number constraints, we have the upper bound $k\leq 8$. However, this is not optimal because there is also the tetrahedral intersection number constraint. This arises because the maximal tetrahedral intersection number for $d=7$ is zero. Taking into account this constraints, we have seen that $k\leq 6$ in this number of dimensions. As before, we consider different values of $k$ case by case. 

\subsection{$k=2$} 
The occupation numbers are $n_2, n_1, n_0$, which satisfy $n_2=n_0$. The only equation is then $2n_2 + n_1= 14$. We also have $n_2=e$, the dimension of the common null subspace. This dimension can take values $7-4=3$ and $7-5=1$. In the first case we have $n_2=n_0=3, n_1=8$. A possible representative of the solution is
\be
\begin{tikzpicture}
\coordinate [label=left:A] (A) at (-1,0);
\coordinate [label=right:B] (B) at (1,0);

\draw[blue, very thick] (A) -- (B);
\node[fill=white] (c) at ($(B)!0.5!(A)$) {$3$};

\end{tikzpicture} \qquad
\raise2ex\vbox{
\hbox{$\psi_A=a_1 a_2 a_3 a_4 a_5 a_6 a_7,$}
\hbox{$\psi_B=a_1^\dagger a_2^\dagger  a_3^\dagger a_4^\dagger a_5 a_6 a_7.$}}
\ee
This is the already familiar ``effectively" 8D solution. The simple part of the stabiliser is ${\rm Spin}(7)\times {\rm SL}(3)$. This is the case $b$ in \cite{Popov-short}.

When the dimension of the common null subspace is $e=1$ we have $n_2=n_0=1$ and $n_1=12$. A possible pair of spinors realising this solution is
\be
\begin{tikzpicture}
\coordinate [label=left:A] (A) at (-1,0);
\coordinate [label=right:B] (B) at (1,0);

\draw[blue, very thick] (A) -- (B);
\node[fill=white] (c) at ($(B)!0.5!(A)$) {$1$};

\end{tikzpicture} \qquad
\raise2ex\vbox{
\hbox{$\psi_A=a_1 a_2 a_3 a_4 a_5 a_6 a_7,$}
\hbox{$\psi_B=a_1^\dagger a_2^\dagger  a_3^\dagger a_4^\dagger a_5^\dagger a_6^\dagger a_7.$}}
\ee
We recognise here the 12D spinor (\ref{12d-p2-2}), where an additional pair of oscillators have been added, but plays no role. The stabiliser of this spinor contains ${\rm SL}(6)$. This is the case $d$ in \cite{Popov-short}.

\subsection{$k=3$}

There is a single independent equation $n_3+ n_2=7$ in this case, 
which in particular implies $n_2\leq 7$. As before, we envisage drawing a triangle with the 3 spinors sitting at its vertices, and putting the pairwise intersection numbers on the edges. It is clear that
\be\label{intersections3}
3 n_3 + n_2 = \sum_{i=1}^3 e_i,
\ee
which, rewritten in terms of $n_3$ becomes
\be
n_3 = \frac{1}{2}(\sum_{i=1}^3 e_i- 7).
\ee 
\subsubsection{The sum of intersection numbers equals 9}
The maximal possible value of each intersection dimension is $3$, and so the maximal value of the right-hand side and thus $n_3$ is $n_3=1$. Then $n_2=6$. A possible triple of pure spinors realising this solution is
\be\label{14d-p3-1}
\begin{tikzpicture}
\coordinate [label=left:A] (A) at (-2*0.75,0);
\coordinate [label=above:B] (B) at (0,3.46*0.75);
\coordinate [label=right:C] (C) at (2*0.75,0);

\draw[blue, very thick] (A) -- (B);
\node[fill=white] (c) at ($(B)!0.5!(A)$) {$3$};

\draw[blue, very thick] (B) -- (C);
\node[fill=white] (c) at ($(C)!0.5!(B)$) {$3$};

\draw[blue, very thick] (A) -- (C);
\node[fill=white] (c) at ($(C)!0.5!(A)$) {$3$};
\end{tikzpicture} \qquad
\raise6ex\vbox{
\hbox{$\psi_A=a_1 a_2 a_3 a_4 a_5 a_6 a_7,$}
\hbox{$\psi_B=a_1^\dagger a_2^\dagger a_3^\dagger a_4^\dagger a_5 a_6 a_7,$}
\hbox{$\psi_C=a_1 a_2  a_3^\dagger a_4^\dagger a_5^\dagger a_6^\dagger a_7.$}}
\ee
It is clear that this is just the 12D solution (\ref{12d-p3}), with one extra fermionic oscillator added. Hence the simple part of the stabiliser is given by ${\rm Sp}(6)$. This is the case $c$ in \cite{Popov-short}.

\subsubsection{The sum of intersection numbers equals 7}
The only other possible solution in the impurity three case is $e_1=e_2=3, e_3=1$, implying the total sum of intersection numbers to be 7. This gives $n_3=n_0=0$, and thus $n_2=n_1=7$. A possible representative of this solution is 
\be\label{14d-p3-2}
\begin{tikzpicture}
\coordinate [label=left:A] (A) at (-2*0.75,0);
\coordinate [label=above:B] (B) at (0,3.46*0.75);
\coordinate [label=right:C] (C) at (2*0.75,0);

\draw[blue, very thick] (A) -- (B);
\node[fill=white] (c) at ($(B)!0.5!(A)$) {$3$};

\draw[blue, very thick] (B) -- (C);
\node[fill=white] (c) at ($(C)!0.5!(B)$) {$3$};

\draw[blue, very thick] (A) -- (C);
\node[fill=white] (c) at ($(C)!0.5!(A)$) {$1$};
\end{tikzpicture} \qquad
\raise6ex\vbox{
\hbox{$\psi_A=a_1 a_2 a_3 a_4 a_5 a_6 a_7,$}
\hbox{$\psi_B=a_1^\dagger a_2^\dagger a_3^\dagger a_4^\dagger a_5 a_6 a_7,$}
\hbox{$\psi_C=a_1^\dagger a_2^\dagger a_3^\dagger  a_4 a_5^\dagger a_6^\dagger a_7^\dagger.$}}
\ee
It is again sufficient to use $e^{\lambda_1 t_1}$ and $e^{\lambda_5 t_5}$ to rescale the coefficients in front of all 3 pure spinors equal. 

To understand the stabiliser we compute
\be
B_{3}(\psi)=\bar{e}_5\wedge \bar{e}_6\wedge \bar{e}_7 + e_1\wedge e_2 \wedge e_3 - \frac{1}{2} \bar{e}_4 \wedge \lambda,
\ee
where $\lambda= \bar{e}_1 \wedge e_1+ \ldots \bar{e}_7\wedge e_7$. The first term in $B_3(\psi)$ is $B_3(\psi_A,\psi_B)$, the second is $B_3(\psi_B,\psi_C)$, and the last is
$B_3(\psi_A,\psi_C)$. It is clear that the last term in $B_3(\psi)$ is stablised by ${\rm SL}(6)$ that acts in the 12D space spanned by $\bar{e}_{1,2,3,5,6,7}$ and $e_{1,2,3,5,6,7}$, thus not touching the directions $\bar{e}_4, e_4$. This ${\rm SL}(6)$ is further broken down to ${\rm SL}(3)\times{\rm SL}(3)$ acting in the spaces spanned by $\bar{e}_{1,2,3}, e_{1,2,3}$ and $\bar{e}_{5,6,7}, e_{5,6,7}$. Hence the simple part of the stabiliser is ${\rm SL}(3)\times{\rm SL}(3)$, reproducing the case $f$ in \cite{Popov-short}.

\subsection{$k=4$}

The equation relating the occupation numbers in this case reads $2n_4 + 2 n_3+ n_2= 14$. The relation to the intersection numbers of the edges is
\be\label{inters-4}
6 n_4 + 3 n_3 + n_2 = \sum_{i=1}^6 e_i.
\ee
Expressing $n_2=14- 2n_3-2n_4$ we get
\be
4n_4+n_3 = \sum_{i=1}^6 e_i-14.
\ee
\subsubsection{Sum of intersection numbers equals 18}
When the sum of intersection numbers is the maximal possible one -- all intersections are 3 -- we have two possible solutions. In the first case $n_4=n_0=0, n_3=n_1=4$, and then $n_2= 6$. A possible representative of this solution is 
\be\label{16d-p4-1}
\begin{tikzpicture}
\coordinate [label=left:A] (A) at (-1.5,0);
\coordinate [label=above:B] (B) at (0,1.5);
\coordinate [label=right:C] (C) at (1.5,0);
\coordinate [label=below:D] (D) at (0,-1.5);

\draw[blue, very thick] (A) -- (B);
\node[fill=white] (c) at ($(B)!0.5!(A)$) {$3$};

\draw[blue, very thick] (B) -- (C);
\node[fill=white] (c) at ($(C)!0.5!(B)$) {$3$};

\draw[blue, very thick] (C) -- (D);
\node[fill=white] (c) at ($(D)!0.5!(C)$) {$3$};

\draw[blue, very thick] (D) -- (A);
\node[fill=white] (c) at ($(A)!0.5!(D)$) {$3$};

\draw[blue, very thick] (A) -- (C);
\node[fill=white] (c) at ($(C)!0.75!(A)$) {$3$};

\draw[blue, very thick] (B) -- (D);
\node[fill=white] (c) at ($(D)!0.75!(B)$) {$3$};
\end{tikzpicture}\qquad
\raise6ex\vbox{
\hbox{$\psi_A=a_1 a_2 a_3 a_4 a_5 a_6 a_7,$ }
\hbox{$\psi_B=a_1^\dagger a_2^\dagger a_3^\dagger a_4^\dagger a_5 a_6  a_7,$}
\hbox{$\psi_C=a_1^\dagger a_2 a_3^\dagger a_4 a_5 a_6^\dagger a_7^\dagger,$}
\hbox{$\psi_D=a_1 a_2^\dagger a_3^\dagger a_4 a_5^\dagger a_6^\dagger  a_7.$}}
\ee
To show that all relative coefficients in this case can be rescaled away, it suffices to consider the action of $e^{\lambda_1 t_1}, e^{\lambda_2 t_2},e^{\lambda_3 t_3}$. The desired transformation is then
\be
c = c_A e^{\lambda_1+\lambda_2+\lambda_3} = c_B e^{-\lambda_1-\lambda_2-\lambda_3}= c_C e^{-\lambda_1+\lambda_2-\lambda_3}=c_De^{\lambda_1-\lambda_2-\lambda_3},
\ee
which is clearly possible. 

To understand the stabiliser we compute $B_{3}(\psi)$. It is given by 6 terms, each coming from a pair of pure spinors, and each given by the product of null directions shared by the pair. There are also some relative signs. We get
\be
B_3(\psi)=\bar{e}_5\wedge \bar{e}_6 \wedge \bar{e}_7 - \bar{e}_2\wedge \bar{e}_4 \wedge \bar{e}_5 + \bar{e}_1\wedge \bar{e}_4 \wedge \bar{e}_7 \\ \nonumber
- e_1 \wedge e_3 \wedge \bar{e}_5 - e_2\wedge e_3 \wedge \bar{e}_7 - e_3 \wedge \bar{e}_4 \wedge e_6 .
\ee
Introducing the notation
\be
v^T=
\begin{pmatrix}
    e_3\wedge \bar{e}_5& e_3 \wedge \bar{e}_7 &e_3 \wedge \bar{e}_4 & \bar{e}_7 \wedge \bar{e}_4 &\bar{e}_4 \wedge \bar{e}_5 & \bar{e}_5 \wedge \bar{e}_7
    \end{pmatrix},\\
\nonumber
w^T=
\begin{pmatrix}
    e_1 &e_2& e_6 & \bar{e}_1& \bar{e}_2& \bar{e}_6
\end{pmatrix}
\ee
we can rewrite
\be
B_3(\psi) = - w^T \wedge v.
\ee
This makes it clear that we have ${\rm SL}(4)$ acting in the null space spanned by $e_3, \bar{e}_4, \bar{e}_5, \bar{e}_7$ that acts on the 2-forms $v$ as an ${\rm SO}(6)$ transformation (preserving the wedge-product metric in the space of 2-forms). This can be followed by an ${\rm SO}(6)$ transformation on the space spanned by $\bar{e}_{1,2,6}, e_{1,2,6}$, leaving the 3-form invariant. This shows that the stabiliser contains ${\rm SL}(4)$, reproducing the case $e$ in \cite{Popov-short}.

In the other case $n_4=1, n_3=0$ and then $n_2=12$. It is clear that in this solution all 4 pure spinors share a common null direction. In other words, the tetrahedral intersection number for a unique tetrahedron in this case is $t=1$. However, this is not allowed by our tetrahedral intersection constraint $t\leq d-7$. Alternatively, we are effectively dealing with a impurity 4 spinor in 12D. As we know, this impurity four spinor is actually a impurity two spinor. And so the solution in this case is in the same ${\rm Spin}(14)$ orbit as the impurity two solution with intersection number one. We won't consider it any further. 

Let us also find solutions with the edge intersection numbers less than maximal. 
\subsubsection{Sum of intersection numbers equals 16}
When there is a single intersection number one we have $\sum_{i=1}^6 e_i = 16$. In this case, we necessarily have $n_4=0$ and then $n_3= 2, n_2= 10$. A possible representative of this solution is 
\be\label{16d-p4-2}
\begin{tikzpicture}
\coordinate [label=left:A] (A) at (-1.5,0);
\coordinate [label=above:B] (B) at (0,1.5);
\coordinate [label=right:C] (C) at (1.5,0);
\coordinate [label=below:D] (D) at (0,-1.5);

\draw[blue, very thick] (A) -- (B);
\node[fill=white] (c) at ($(B)!0.5!(A)$) {$3$};

\draw[blue, very thick] (B) -- (C);
\node[fill=white] (c) at ($(C)!0.5!(B)$) {$3$};

\draw[blue, very thick] (C) -- (D);
\node[fill=white] (c) at ($(D)!0.5!(C)$) {$3$};

\draw[blue, very thick] (D) -- (A);
\node[fill=white] (c) at ($(A)!0.5!(D)$) {$1$};

\draw[blue, very thick] (A) -- (C);
\node[fill=white] (c) at ($(C)!0.75!(A)$) {$3$};

\draw[blue, very thick] (B) -- (D);
\node[fill=white] (c) at ($(D)!0.75!(B)$) {$3$};
\end{tikzpicture}\qquad
\raise6ex\vbox{
\hbox{$\psi_A=a_1 a_2 a_3 a_4 a_5 a_6 a_7,$ }
\hbox{$\psi_B=a_1^\dagger a_2^\dagger a_3^\dagger a_4 a_5 a_6 a_7^\dagger,$}
\hbox{$\psi_C=a_1 a_2 a_3^\dagger a_4^\dagger a_5^\dagger a_6 a_7^\dagger,$}
\hbox{$\psi_D=a_1^\dagger a_2^\dagger a_3^\dagger  a_4^\dagger a_5^\dagger a_6^\dagger a_7.$}}
\ee

To prove that the relative coefficients in the linear combination of these pure spinors can be rescaled away, we consider the action of $e^{\lambda_2 t_2}, e^{\lambda_3 t_3}, e^{\lambda_4 t_4}$. We want to set all the coefficients to be equal by these transformations, which gives the equations
\be
c= c_A e^{\lambda_2+\lambda_3+\lambda_4} = c_B e^{-\lambda_2-\lambda_3+\lambda_4}=c_C e^{\lambda_2-\lambda_3-\lambda_4}= c_D e^{-\lambda_2-\lambda_3-\lambda_4}.
\ee
This has a solution, and so it is sufficient to consider just the sum of the pure spinors $\psi_{A,B,C,D}$. 

As before, the 3-form $B_{3}(\psi)$ is given by the sum of 6 terms. Five of them are decomposable 3-forms each given by the product of the 3 null directions shared by a pair of pure spinors. In the case of the pure spinors $\psi_A, \psi_D$ that only share a single null direction, the 3-form $B_3(\psi_A, \psi_D)$ is given by the product of the shared null direction, times the 2-form for the remaining 12D space
\be
B_3(\psi)=(\bar{e}_1\wedge \bar{e}_2 + \bar{e}_4\wedge \bar{e}_5)\wedge \bar{e}_6 - (e_1\wedge e_2 + e_4 \wedge e_5)\wedge e_6
\\ \nonumber
- e_3 \wedge \bar{e}_6 \wedge e_7 + \bar{e}_7\wedge \lambda.
\ee
It is clear that the stabiliser contains the group ${\rm Sp}_4$ acting on the space spanned by $e_{1,2,4,5}$ and preserving the 2-form $e_1\wedge e_2 + e_4 \wedge e_5$, and similarly in the space spanned by $\bar{e}_{1,2,4,5}$ and preserving $\bar{e}_1\wedge \bar{e}_2 + \bar{e}_4\wedge \bar{e}_5$. There is also a copy of ${\rm SL}_2$ in the stabiliser, which is not easy to see geometrically. However, the stabiliser subalgebra is easily computed explicitly. One obtains that the stabiliser contains $\chi_{{\rm SL}_2}+\chi_{{\rm Sp}_4}$, where these subalgebras are given by
\be
\begin{split}
    \chi_{{\rm Sp}_4}=&\{a_1a_4^{\dagger}-a_5a_2^{\dagger}, a_4a_1^{\dagger}-a_2a_5^{\dagger} ,a_5a_1^{\dagger}+a_2a_4^{\dagger}, a_1a_5^{\dagger}+a_4a_2^{\dagger}, a_1a_1^{\dagger}-a_2a_2^{\dagger}, a_4a_4^{\dagger}-a_5a_5^{\dagger},\\ 
    &a_1a_2^{\dagger}, a_2a_1^{\dagger}, a_4a_5^{\dagger}, a_5a_4^{\dagger}\}, \ {\rm and}\ \\
    \chi_{{\rm SL}_2}=&\{a_1a_1^{\dagger}+a_2a_2^{\dagger}+a_4a_4^{\dagger}+a_5a_5^{\dagger}-2a_3a_3^{\dagger}-2a_6a_6^{\dagger},2a_6a_7^{\dagger}-a_7a_3+a_1^{\dagger}a_2^{\dagger}+a_4^{\dagger}a_5^{\dagger},\\
    &-a_7a_6^{\dagger}+2a_3^{\dagger}a_7^{\dagger}+a_1a_2+a_4a_5\}.
\end{split}
\ee
It can easily be shown that these two subalgebras commute. Thus the stabilising subgroup contains ${\rm SL}_2\times{\rm Sp}_4$, reproducing the case $g$ in \cite{Popov-short}.

\subsubsection{Sum of intersection numbers equals 14}
When there are two intersection numbers one then $n_4=n_3=0, n_2=14$. A possible representative of this solution is
\be\label{16d-p4-3}
\begin{tikzpicture}
\coordinate [label=left:A] (A) at (-1.5,0);
\coordinate [label=above:B] (B) at (0,1.5);
\coordinate [label=right:C] (C) at (1.5,0);
\coordinate [label=below:D] (D) at (0,-1.5);

\draw[blue, very thick] (A) -- (B);
\node[fill=white] (c) at ($(B)!0.5!(A)$) {$3$};

\draw[blue, very thick] (B) -- (C);
\node[fill=white] (c) at ($(C)!0.5!(B)$) {$1$};

\draw[blue, very thick] (C) -- (D);
\node[fill=white] (c) at ($(D)!0.5!(C)$) {$3$};

\draw[blue, very thick] (D) -- (A);
\node[fill=white] (c) at ($(A)!0.5!(D)$) {$1$};

\draw[blue, very thick] (A) -- (C);
\node[fill=white] (c) at ($(C)!0.75!(A)$) {$3$};

\draw[blue, very thick] (B) -- (D);
\node[fill=white] (c) at ($(D)!0.75!(B)$) {$3$};
\end{tikzpicture}\qquad
\raise6ex\vbox{
\hbox{$\psi_A=a_1 a_2 a_3 a_4 a_5 a_6 a_7,$ }
\hbox{$\psi_B=a_1^\dagger a_2^\dagger a_3^\dagger a_4^\dagger a_5 a_6 a_7,$}
\hbox{$\psi_C=a_1 a_2 a_3 a_4^\dagger a_5^\dagger a_6^\dagger a_7^\dagger,$}
\hbox{$\psi_D=a_1^\dagger a_2^\dagger a_3^\dagger  a_4 a_5^\dagger a_6^\dagger a_7^\dagger.$}}
\ee

We again show that the relative coefficients are irrelevant, by considering the action of the Cartan generators $e^{\lambda_3 t_3}, e^{\lambda_4 t_4}, e^{\lambda_5 t_5}$. The equations setting transformed coefficients equal to each other are
\be
c= c_A e^{\lambda_3+\lambda_4+\lambda_5} = c_B e^{-\lambda_3-\lambda_4+\lambda_5}=c_C e^{\lambda_3-\lambda_4-\lambda_5}=c_D e^{-\lambda_3+\lambda_4-\lambda_5}.
\ee
These equations have a solution, and so the relative coefficients can be scaled away.  

The 3-form is computed by considering pairwise intersections of the pure spinors involved. We get
\be
B_3(\psi)=- \frac{1}{2}\bar{e}^4 \wedge \lambda+\bar{e}^1 \wedge \bar{e}^2 \wedge \bar{e}^3 - e^1 \wedge e^2 \wedge e^3 \\
\nonumber
+\frac{1}{2}e^4\wedge\bar{\lambda} +\bar{e}^5 \wedge \bar{e}^6 \wedge \bar{e}^7 + e^5 \wedge e^6 \wedge e^7 ,
\ee
where
\be
\lambda := \bar{e}^1\wedge e^1 + \bar{e}^2\wedge e^2+\bar{e}^3\wedge e^3+ \bar{e}^5 \wedge e^5 + \bar{e}^6\wedge e^6+\bar{e}^7\wedge e^7, \\ \nonumber
\bar{\lambda}:=-\bar{e}^1\wedge e^1 - \bar{e}^2\wedge e^2-\bar{e}^3\wedge e^3+ \bar{e}^5 \wedge e^5 + \bar{e}^6\wedge e^6+\bar{e}^7\wedge e^7.
\ee
This can be rewritten as
\be
B_3(\psi)=- \frac{1}{2}(e^4+\bar{e}^4) \wedge(\bar{e}^1\wedge e^1 + \bar{e}^2\wedge e^2+\bar{e}^3\wedge e^3)+\bar{e}^1 \wedge \bar{e}^2 \wedge \bar{e}^3 - e^1 \wedge e^2 \wedge e^3 \\
\nonumber
+\frac{1}{2}(e^4-\bar{e}^4)\wedge(\bar{e}^5 \wedge e^5 + \bar{e}^6\wedge e^6+\bar{e}^7\wedge e^7) +\bar{e}^5 \wedge \bar{e}^6 \wedge \bar{e}^7 + e^5 \wedge e^6 \wedge e^7 .
\ee
This is a sum of two canonical ${\rm G}_2$-invariant 3-forms, one in the space spanned by $e^4+\bar{e}^4, e^{1,2,3}, \bar{e}^{1,2,3}$, the other in $e^4-\bar{e}^4, e^{5,6,6}, \bar{e}^{5,6,7}$. Thus, the stabiliser contains ${\rm G}_2\times{\rm G}_2$. This is the generic orbit that is dense in the space of semi-spinors. The degree eight invariant that exists in the space of semi-spinors is different from zero on this orbit. This is the case $i$ in \cite{Popov-short}.

\subsection{$k=5$}

The relation between the occupation numbers is $n_5+n_4+n_3=7$. The relation with the edge intersection numbers is
\be
10 n_5 + 6 n_4 + 3n_3 + n_2 = \sum_{i=1}^{10} e_i,
\ee
where of course $n_2=n_3$. We eliminate $n_3=7-n_5-n_4$ to get
\be
6n_5 + 2n_4 =\sum_{i=1}^{10} e_i -  28.
\ee
This immediately shows that the sum of the intersection numbers on the edges, the first term on the right-hand side, must be at least 28 to have a solution. 

The largest possible value of the sum of intersection numbers is 30. This means that necessarily $n_5=0$. We then have a possible solution with all 10 intersection numbers equal to 3 and $n_4=1$. However, $n_4=1$ means that there are four pure spinors that all share a single common null direction. As we know, this means that these four pure spinors effectively live in 12D, where they correspond to a reducible configuration of impurity two. So, this is not a new orbit. 

The other possibility is when there are 9 edges with intersection number 3, and the last remaining one with intersection number one. This gives the solution $n_5=n_4=n_0=n_1=0, n_3=n_2=7$. A possible representative of this solution is
\be\label{14d-p5}
\begin{tikzpicture}
\coordinate [label=left:A] (A) at (-1.5,0);
\coordinate [label=above:B] (B) at (-0.46,1.42);
\coordinate [label=right:C] (C) at (1.21,0.88);
\coordinate [label=below:D] (D) at (1.21,-0.88);
\coordinate [label=below:E] (E) at (-0.46,-1.42);

\draw[blue, very thick] (A) -- (B);
\node[fill=white] (c) at ($(B)!0.5!(A)$) {$3$};

\draw[blue, very thick] (B) -- (C);
\node[fill=white] (c) at ($(C)!0.5!(B)$) {$3$};

\draw[blue, very thick] (C) -- (D);
\node[fill=white] (c) at ($(D)!0.5!(C)$) {$3$};

\draw[blue, very thick] (D) -- (A);
\node[fill=white] (c) at ($(A)!0.5!(D)$) {$3$};

\draw[blue, very thick] (A) -- (C);
\node[fill=white] (c) at ($(C)!0.5!(A)$) {$3$};

\draw[blue, very thick] (B) -- (D);
\node[fill=white] (c) at ($(D)!0.5!(B)$) {$3$};

\draw[blue, very thick] (A) -- (E);
\node[fill=white] (c) at ($(A)!0.5!(E)$) {$1$};

\draw[blue, very thick] (B) -- (E);
\node[fill=white] (c) at ($(B)!0.5!(E)$) {$3$};

\draw[blue, very thick] (C) -- (E);
\node[fill=white] (c) at ($(C)!0.5!(E)$) {$3$};

\draw[blue, very thick] (D) -- (E);
\node[fill=white] (c) at ($(D)!0.5!(E)$) {$3$};

\end{tikzpicture}\qquad
\raise6ex\vbox{
\hbox{$\psi_A=a_1 a_2 a_3 a_4 a_5 a_6 a_7,$ }
\hbox{$\psi_B=a_1^\dagger a_2^\dagger a_3 a_4^\dagger a_5 a_6 a_7^\dagger,$}
\hbox{$\psi_C=a_1^\dagger a_2 a_3^\dagger a_4^\dagger a_5 a_6^\dagger a_7 ,$}
\hbox{$\psi_D=a_1 a_2^\dagger a_3^\dagger  a_4^\dagger a_5^\dagger a_6 a_7,$}
\hbox{$\psi_E=a_1^\dagger a_2^\dagger a_3^\dagger  a_4 a_5^\dagger a_6^\dagger a_7^\dagger.$}}
\ee
We can again show that all relative coefficients can be rescaled away. To do this, we consider the action of the Cartan generators $e^{\lambda_1 t_1}, e^{\lambda_2 t_2},e^{\lambda_3 t_3},e^{\lambda_4 t_4}$. The desired conditions are
\be\nonumber
c= c_A e^{\lambda_1+\lambda_2+\lambda_3+\lambda_4} = c_B e^{-\lambda_1-\lambda_2+\lambda_3-\lambda_4}=c_C e^{-\lambda_1+\lambda_2-\lambda_3-\lambda_4}=c_D e^{\lambda_1-\lambda_2-\lambda_3-\lambda_4}=c_E e^{-\lambda_1-\lambda_2-\lambda_3+\lambda_4}.
\ee
There is clearly a solution, and so we only need to consider the sum of these pure spinors. It is convenient to take the spinor $\psi_C$ in the form $\psi_C= e_3 \wedge e_1 \wedge e_4 \wedge e_6$, while all other spinors are given with the same signs as one obtains by acting with the creation operators in the order indicated. E.g. $\psi_B = e_1\wedge e_2\wedge e_4\wedge e_7$. A computation shows that $B_{3}(\psi)$ is given as
\be
\begin{split}
    B_3(\psi)&=-\frac{1}{2}\bar{e}_4\wedge \lambda+\bar{e}_1\wedge\bar{e}_6\wedge\bar{e}_7+\bar{e}_2\wedge\bar{e}_7\wedge\bar{e}_5+\bar{e}_3\wedge\bar{e}_5\wedge\bar{e}_6 \\ 
    &\quad \qquad +e_2\wedge e_3 \wedge e_5+ e_3\wedge e_1 \wedge e_6+e_1\wedge e_2 \wedge e_7 \\
    &\quad \qquad+e_4\wedge(\bar{e}_5\wedge e_1 + \bar{e}_6 \wedge e_2 +\bar{e}_7\wedge e_3 ),
\end{split}
\ee
where $\lambda = \bar{e}_1\wedge e_1 +\ldots \bar{e}_7\wedge e_7$. We can now introduce
\be
n_1 = e_1+e_5, \quad n_5 = e_1 - e_5, \quad \bar{n}_1= \bar{e}_1+\bar{e}_5, \quad \bar{n}_5= \bar{e}_1-\bar{e}_5, \\ \nonumber
n_2 = e_2+e_6, \quad n_6 = e_2 - e_6, \quad \bar{n}_2= \bar{e}_2+\bar{e}_6, \quad \bar{n}_6= \bar{e}_2-\bar{e}_6, \\ \nonumber
n_3 = e_3+e_7, \quad n_7 = e_3 - e_7, \quad \bar{n}_3= \bar{e}_3+\bar{e}_7, \quad \bar{n}_7= \bar{e}_3-\bar{e}_7.
\ee
In terms of these new basis the 3-form becomes
 \be
 B_3(\psi)=\frac{1}{2}(e_4-\bar{e}_4)\wedge ( \bar{n}_1 \wedge n_1 +\bar{n}_2\wedge n_2 + \bar{n}_3\wedge n_3)+ \frac{1}{2}n_1\wedge n_2\wedge n_3+ \frac{1}{2}\bar{n}_1\wedge \bar{n}_2\wedge \bar{n}_3
 \\ \nonumber
 - \frac{1}{2} (e_4+\bar{e}_4) \wedge(  \bar{n}_5\wedge n_5 +\bar{n}_6\wedge n_6 +\bar{n}_7\wedge n_7 )
  - \frac{1}{2}n_5\wedge n_6\wedge n_7  + \frac{1}{2}\bar{n}_5\wedge \bar{n}_6\wedge \bar{n}_7  \\ \nonumber
- \frac{1}{8}(n_1-n_5)\wedge (n_2-n_6)\wedge (n_3-n_7) - \frac{1}{8}(\bar{n}_1+\bar{n}_5)\wedge (\bar{n}_2+\bar{n}_6)\wedge (\bar{n}_3+\bar{n}_7).
 \ee
 The first two lines are the canonical 3-forms in the spaces spanned by $e_4-\bar{e}_4, n_{1,2,3}, \bar{n}_{1,2,3}$ and $e_4+\bar{e}_4, n_{5,6,7}, \bar{n}_{5,6,7}$. Each of them is invariant under a copy of $G_2$ acting on the corresponding space. The last line breaks this pair of $G_2$ groups to just a single copy of $G_2$ acting simultaneously on all the directions. We thus reproduce the orbit that contains $G_2$ in the stabiliser. This is the case $h$ in \cite{Popov-short}. This is the orbit on which the degree eight invariant existing in the space of semi-spinors vanishes, and which contains all other orbits (except the generic one (\ref{16d-p4-3})) as its degenerations. Indeed, removing the edge labelled with $1$, together with one of the two vertices sharing this edge, one gets the orbit (\ref{16d-p4-1}). Removing instead one of the edges labelled $3$, together with one of its vertices, one gets the orbit (\ref{16d-p4-2}). It is clear from the diagram that the generic orbit (\ref{16d-p4-3}) cannot be obtained as a degeneration of (\ref{14d-p5}). 

\subsection{$k=6$}

The relation between the occupation numbers gives $n_3 = 14 - 2n_6 - 2n_5 - 2n_4$.
The relation with the intersection numbers is
\be
15 n_6+ 10 n_5 + 7 n_4 + 3n_3 = \sum_{i=1}^{15} e_i,
\ee
where $n_4=n_2$. Substituting the expression for $n_3$ in terms of $n_{4,5,6}$ we get
\be
9 n_6+ 4 n_5 +  n_4  = \sum_{i=1}^{15} e_i - 42.
\ee
The largest possible value of the sum of the intersection numbers is 45. This shows that $n_5=n_6=0$. When the sum of the intersection numbers is 45 we get $n_4=n_2=3$ and $n_3=8$. However, as we already know from the previous considerations, we cannot have $n_4$ different from zero, as then there are tetrahedra with tetrahedral intersection numbers different from zero. So, the sum of edge intersection numbers must be 42. However, this is not possible, because its possible values are $45, 43, 41, \ldots$. So, there is no solution with $k=6$ that gives rise to an irreducible configuration of pure spinors in this case. 

We do not need to consider configurations with large values of $k$ because these, even though possible by edge intersection constraints, necessarily have tetrahedral intersections different from zero, or equivalently $n_4\not=0$, as can be easily checked. They correspond to reducible configurations of pure spinors, in the sense that an impure spinor of impurity $k+2$ is in the same orbit as an impure spinor of impurity $k$. We do not need to consider these cases as no new spinor orbits arise this way.

\end{document}